\documentclass[12pt, a4paper, reqno]{amsart}
\usepackage{amscd, amssymb,amsmath}
\allowdisplaybreaks[3]
\begin{document}

\def\supp{\operatorname{supp}}
\def\Ind{\operatorname{Ind}}
\def\Aut{\operatorname{Aut}}
\def\max{\operatorname{max}}
\def\Rep{\operatorname{Rep}}
\def\sp{\operatorname{sp}}
\def\clsp{\overline{\operatorname{sp}}}
\def\dashind{\operatorname{\!-Ind}}
\def\id{\operatorname{id}}
\def\m{\operatorname{m}}
\def\rt{\operatorname{rt}}
\def\lt{\operatorname{lt}}
\def\H{\mathcal{H}}
\def\K{\mathcal{K}}
\def\C{\mathbb{C}}
\def\T{\mathbb{T}}
\def\justMN{M}
\def\justMH{M}

\newcommand{\kaz}[1]{{#1}}

\newtheorem{thm}{Theorem}  
\newtheorem{cor}[thm]{Corollary}
\newtheorem{lemma}[thm]{Lemma}
\newtheorem{prop}[thm]{Proposition}
\newtheorem{thm1}{Theorem}

\theoremstyle{definition}
\newtheorem{defn}[thm]{Definition}
\newtheorem{remark}[thm]{Remark}
\newtheorem{example}[thm]{Example}
\newtheorem{remarks}[thm]{Remarks}
\newtheorem{claim}[thm]{Claim}
\newtheorem{problem}[thm]{Problem}
\newtheorem*{problem*}{Problem}

\numberwithin{equation}{section}

\title[Extension problems and non-abelian duality for $C^*$-algebras]{\boldmath Extension problems and
\\non-abelian duality for $C^*$-algebras}

\author[an Huef]{Astrid an Huef}
\address{School of Mathematics and Statistics\\
The University of New South Wales\\
NSW 2052\\
Australia}
\email{astrid@unsw.edu.au}

\author[Kaliszewski]{S. Kaliszewski}
\address{Department of Mathematics 
\kaz{and Statistics}
\\Arizona State University\\ AZ
85287-1804\\USA} \email{kaliszewski@asu.edu}

\author[Raeburn]{Iain Raeburn}
\address{School of Mathematical and Physical Sciences\\
University of Newcastle\\
NSW 2308\\ Australia}
\email{iain.raeburn@newcastle.edu.au}

\thanks{This research was supported by grants from the Australian Research Council, the National Science Foundation and the University of New South Wales.}


\date{September 8, 2006}

\begin{abstract}
Suppose that $H$ is a closed subgroup of a locally compact group $G$. We show that a unitary representation $U$ of $H$ is the restriction of a unitary representation of $G$ if and only if a dual representation $\widehat U$ of a crossed product $C^*(G)\rtimes (G/H)$ is regular in an appropriate sense. We then discuss the problem of deciding whether a given representation is regular; we believe that this problem will prove to be an interesting test question in non-abelian duality for crossed products of $C^*$-algebras.
\end{abstract}

\maketitle

In \cite{aHRmackey}, we considered the classical problem of deciding
when a unitary representation $U$ of a closed normal subgroup $H$ of a
locally compact group $G$ 
\kaz{extends to a} 
unitary representation of
$G$ on the same space. We proved that $U$ extends to $G$ if and only
if a dual representation $\widehat U$ of a coaction crossed product
$C^*(G)\rtimes(G/H)$ is regular, and used Mansfield's imprimitivity
theorem \cite{man} to describe these regular representations when $G$ is
amenable \cite[Theorem~5]{aHRmackey}.  
\kaz{Since regular representations
are induced representations, this accords with} 
a general principle in
non-abelian duality for crossed products of $C^*$-algebras:
 duality swaps restriction and induction of representations
\cite{E,KQR,BE}.
The extension problem makes sense for an
arbitrary (not necessarily normal) closed subgroup $H$, and finding
the appropriate dual formulation is an intriguing test question for any
theory of crossed products by homogeneous spaces.

In general we do not know  what it means to say that a homogeneous 
\kaz{space~$G/H$}
coacts on a $C^*$-algebra, but when we have a coaction of $G$
on a $C^*$-algebra $A$, we can define 
\kaz{a crossed product $C^*$-algebra}
$A\rtimes (G/H)$, and 
\kaz{such}
algebras arise in applications (see,
for example, \cite{EQ,DPR, PQR}). The coaction arising in the above
extension problem is the restriction of a dual coaction of $G$, and the
crossed product $C^*(G)\rtimes (G/H)$ has a natural family 
\kaz{of}
regular
representations. Our Theorem~\ref{thm-ext-reg} says that $U$ extends to
a unitary representation of $G$ if and only if $\widehat U$ is regular
in this sense. The coactions are only needed to motivate the result, which
 can be formulated using only ordinary crossed products by actions,
and its proof uses a direct spatial argument.

Theorem~\ref{thm-ext-reg} raises the question of how one can decide
whether a given representation of a crossed product $C^*(G)\rtimes(G/H)$
is regular.  Our Theorem~\ref{thm-fm} answers this question when $H$ is
normal, and gives a version of \cite[Theorem~5]{aHRmackey} which works
for non-amenable groups (see Corollary~\ref{thm-main}). We have also been
able to make some progress on this question in the purely algebraic
situation where $H$ is not normal but $(G,H)$ is a discrete Hecke pair;
since this has required completely different methods, we  discuss
these results elsewhere \cite{aHKRHecke}.

Our results sit naturally in the context of $C^*$-dynamical systems consisting of an action $\alpha$ of a locally compact group $G$ by automorphisms of a $C^*$-algebra $A$, and we work in this generality throughout.

\section{Preliminaries}

Throughout, $H$ is a closed subgroup of a locally compact group
$G$, and $N$ 
\kaz{is}
a closed normal subgroup.  The group $G$ acts on the
left of the homogeneous space $G/H$, and this induces a continuous
action $\lt:G\to\Aut C_0(G/H)$ defined by $\lt_t(f)(sH)=f(t^{-1}sH)$.
A representation $V$ of $G$ is a continuous homomorphism $V$ of $G$
into the group $U(\H)$ of unitary operators on a Hilbert space $\H$,
where $U(\H)$ has the strong operator topology.

We use left Haar measures, and denote by $\Delta_G$ and $\Delta_H$ the modular functions on $G$ and $H$.  Let $\rho:G\to(0,\infty)$ be a continuous function
such that $\rho(st)=\frac{\Delta_H(t)}{\Delta_G(t)}\rho(s)$ for $s\in G$ and
$t\in H$.  From \cite[\S C.1]{tfb}, for example, there exists a quasi-invariant measure on the quotient space $G/H$ such
that
\[
\int_G f(s)\rho(s)\, ds=\int_{G/H}\int_H f(st)\, dt\, d(sH)\quad\text{for\
}f\in C_c(G). \]
If $N$ is a closed normal subgroup of $G$ with $N\subset H$ then  $\Delta_H(n)=\Delta_N(n)=\Delta_G(n)$ and hence \ $\rho(sn)=\rho(s)$ for
all $n\in N$ and $s\in G$. We can thus choose left Haar measures on $G/N$ and $H/N$ such that
\[
\int_{G/N} f(sN)\rho(s)\, d(sN)=\int_{G/H}\int_{H/N} f(stN)\, d(tN)\,
d(sH)\quad\text{for\ }f\in C_c(G/N).
\]
We denote by
$\lambda^{G/H}$ the quasi-regular representation of $G$ on $L^2(G/H)$ given by
\[
(\lambda^{G/H}_r\xi)(sH)=\sqrt{\frac{\rho(r^{-1}s)}{\rho(s)}}\xi(r^{-1}sH)
\quad \text{for $s,r\in G$ and $\xi\in L^2(G/H)$},
\]
and 
\kaz{we denote}
by $\justMH$ the representation of $C_0(G/H)$ by multiplication operators
on $L^2(G/H)$. We  use  $\lambda$ and $\nu$ to denote  left- and
right-regular representations, respectively; so, for example,  $\lambda$
could denote  the left regular representation of $G$ on 
\kaz{$L^2(G)$, or of $G/N$}
on $L^2(G/N)$.

Our conventions regarding actions and coactions of $G$ on a $C^*$-algebra $A$ are those of \cite[Appendix~A]{BE}. We emphasise that all the group algebras and crossed products in this paper are full; in other words, they are universal for an appropriate family of covariant representations.

We denote by $X$ or $X_H^G(A)$ the $((A\otimes C_0(G/H))\rtimes_{\alpha\otimes\lt}
G)$--$(A\rtimes_\alpha H)$ imprimitivity bimodule implementing Green's imprimitivity theorem \cite[Theorem~6]{green}.  We  view $X$ as the completion of $X_0:=C_c(G,A)$ with respect to the formulas for the actions and inner products from \cite[Equations~(B.5)]{BE}.  Tensoring with $X$ then gives a bijection
\[
X\dashind:\Rep (A\rtimes_\alpha H)\to \Rep((A\otimes C_0(G/H))\rtimes_{\alpha\otimes\lt}G)
\]
on unitary equivalence classes of representations \cite[Theorem~3.29]{tfb}.

\section{Coaction-regular representations and the extension problem}
Suppose $(A,G,\alpha)$ is a dynamical system, $H$ is a closed subgroup of $G$, and $(\pi,U)$ is a covariant representation of $(A,H,\alpha)$. We want to find a condition on the induced representation $X\dashind (\pi\rtimes U)$ of $(A\otimes C_0(G/H))\rtimes_{\alpha\otimes\lt} G$ which is equivalent to saying that $(\pi,U)$ extends to a covariant representation of $(A,G,\alpha)$. Our result will say that $(\pi,U)$ extends if and only if $X\dashind (\pi\rtimes U)$ is a kind of regular representation. However, to see why we call these representations regular we have to view $(A\otimes C_0(G/H))\rtimes_{\alpha\otimes\lt} G$ as a crossed product of $A\rtimes_\alpha G$ by a coaction of the homogeneous space $G/H$. As motivation for this, we look at what happens when we have a closed normal subgroup $N$.

Denote by $\widehat\alpha:A\to M((A\rtimes_\alpha G)\otimes C^*(G))$ the dual coaction of $G$ on the crossed product $A\rtimes_\alpha G$ \cite[A.26]{BE}. Since the subgroup $N$ is normal, we can restrict $\widehat\alpha$ to a coaction $\widehat\alpha|$ of $G/N$  on $A\rtimes_\alpha G$ \cite[A.28]{BE}.
It is proved in \cite[Proposition~A.63 and Theorem~A.64]{BE}, for example, that there is an isomorphism of $(A\rtimes_\alpha G)\rtimes_{\widehat\alpha|}G/N$ onto $(A\otimes C_0(G/N))\rtimes_{\alpha\otimes\lt} G$
which carries 
\kaz{a}
representation $(\phi\rtimes W)\rtimes\mu$ of
$(A\rtimes_\alpha G)\rtimes_{\widehat\alpha|}G/N$ into the
representation $(\phi\otimes\mu)\rtimes W$  of $(A\otimes
C_0(G/N))\rtimes_{\alpha\otimes\lt} G$.  The \emph{regular
representation of $(A\rtimes_\alpha G)\rtimes_{\widehat\alpha|}G/N$
induced from} 
\kaz{a covariant representation $(\phi, W)$ of $(A,G,\alpha)$}
is the representation \[(((\phi\rtimes
W)\otimes\lambda)\circ\widehat\alpha|)\rtimes (1\otimes \justMN)\]
on $\H\otimes L^2(G/N)$, and we call the corresponding representation
\[(\phi\otimes \justMN)\rtimes (W\otimes \lambda^{G/N})\] of $(A\otimes
C_0(G/N))\rtimes_{\alpha\otimes\lt} G$ the \emph{coaction-regular
representation induced from $(\phi, W)$}.

When $H$ is not normal, we do not have a satisfactory notion of coaction of $G/H$, and hence no completely satisfactory notion of crossed product by $G/H$. However, here we have a coaction $\widehat \alpha$ of $G$, and $(A\otimes C_0(G/H))\rtimes_{\alpha\otimes\lt} G$ is one of the candidates for what we should mean by the crossed product of $A\rtimes_\alpha G$ by the coaction $\widehat\alpha$ of the homogeneous space $G/H$ (see the discussion in \cite[\S2]{EKR}). Further, we can still use the quasi-regular representation $\lambda^{G/H}$ of $G$ to build a representation
\[
(\phi\otimes
\justMH)\rtimes (W\otimes \lambda^{G/H}):(A\otimes
C_0(G/H))\rtimes_{\alpha\otimes\lt} G\to B(\H\otimes L^2(G/H)),
\]
\kaz{and, as in the normal case,}
we call this the \emph{coaction-regular representation  induced
from  $(\phi,W)$}. We do this to emphasise that it is only regular
when we think of $(A\otimes C_0(G/H))\rtimes_{\alpha\otimes\lt} G$ as
a crossed product by $G/H$;  the regular representations of $(A\otimes
C_0(G/H))\rtimes_{\alpha\otimes\lt} G$ induced from representations of
$A\otimes C_0(G/H)$ are quite different, and it is the coaction-regular
representations which arise naturally in our dual problem.

\begin{thm}\label{thm-ext-reg} 
Suppose that $\alpha$ is an action of a locally compact group $G$ on
a $C^*$-algebra $A$ and $H$ is a closed subgroup of $G$.  Let $X$ be
the imprimitivity bimodule implementing Green's imprimitivity theorem,
and 
\kaz{let $(\pi, U)$ be}
a covariant representation of $(A,H,\alpha)$ on 
\kaz{a Hilbert space~$\H$.}
Then there exists a covariant representation $(\pi, V)$ of $(A,G,\alpha)$
on $\H$ such that $U=V|_H$ if and only if $X\dashind(\pi\rtimes  U)$
is unitarily equivalent to a coaction-regular representation.
\end{thm}

\begin{proof}
We will show that 
\kaz{for any covariant representation $(\pi,V)$ of $(A,G,\alpha)$,
we have}
\begin{equation}\label{eq-onlyif}
X\dashind(\pi\rtimes V|_H) \sim \big(\pi\otimes \justMH  \big)\rtimes
\big( V\otimes\lambda^{G/H} \big).
\end{equation}
This immediately gives the ``only if'' direction of the theorem. It also gives  the converse: if
$X\dashind(\pi\rtimes U)$ is unitarily equivalent to the coaction-regular
representation induced from 
\kaz{some covariant representation $(\phi,W)$,}
then
\eqref{eq-onlyif} gives
\[
X\dashind(\phi\rtimes W|_H)
\sim
\big(\phi\otimes \justMH \big)\rtimes \big( W\otimes\lambda^{G/H} \big)
\sim
X\dashind(\pi\rtimes U),
\] so that $\pi\rtimes U$ is unitarily equivalent to
$\phi\rtimes W|_H$  because $X$ is an
imprimitivity bimodule.  Now moving $W$ over to the space of $(\pi, U)$ gives
\kaz{the desired representation~$V$.}

So suppose 
\kaz{that $(\pi,V)$ is a representation}
of $(A,G,\alpha)$ 
\kaz{and $U=V|_H$.}
We define
$\Psi: X_0\odot_{A\rtimes_\alpha H}\H\to L^2(G/H,\H)$ by
\begin{equation}\label{defn-Psi}
\Psi(x\otimes h)(sH)=\int_H\pi(x(st))V_{st}h
\Delta_G(st)^{-1/2}\rho(st)^{-1/2}\, dt;
\end{equation}
the usual arguments show that $\Psi(x\otimes h)$ is a continuous function of compact support, and hence belongs to $L^2(G/H,\H)$.
 We claim that $\Psi$ extends to a unitary isomorphism of $X\otimes_{A\rtimes_\alpha H}\H$ onto $L^2(G/H,\H)$ which
intertwines $X\dashind(\pi\rtimes U)$ and  $(\pi\otimes \justMH )\rtimes ( V\otimes\lambda^{G/H} )$.

We start by checking  that $\Psi$ is isometric. Let $x,y\in X_0$ and
$h,k\in\H$. Then using the formulas for the inner product in \cite[(B.5)]{BE}, we have
\begin{align}
(x\otimes h\, |\, y\otimes k)
&=(\pi\rtimes U(\langle y\,, x\, \rangle_{A\rtimes_\alpha
H})h\, |\, k)\notag\\
&=\Big(\Big(\int_H \pi( \langle y\,, x\, \rangle_{A\rtimes_\alpha
H} (t))U_t\, dt\Big) h\, \Big|\, k\Big)\notag\\
&=\int_H\int_G \big(\pi(\alpha_s(y(s^{-1})^*x(s^{-1}t)))U_th\, \big|\,
k\big)\Delta_H(t)^{-1/2}\, ds\, dt\notag\\
&=\int_H\int_G\big(\pi(\alpha_s(x(s^{-1}t)))U_th\,
\big|\,\pi(\alpha_s(y(s^{-1})))k\big) \Delta_H(t)^{-1/2}\, ds\, dt \notag\\
&=\int_H\int_G\big(\pi(x(s^{-1}t))V_{s^{-1}t}h\, \big|\,
\pi(y(s^{-1}))V_{s^{-1}}k  \big) \Delta_H(t)^{-1/2}\, ds\,
dt\label{eq-iso1}
\end{align}
using the covariance of $(\pi,V)$ for $(A,G,\alpha)$ and $U_t=V_t$ for $t\in H$.
At the expense of a modular function we replace $s$ by $s^{-1}$, and then we write the integral over
$G$ as an iterated integral over $G/H$ and $H$. This gives
\begin{align}
\eqref{eq-iso1}
&=\int_H\int_{G/H}\int_H \big(\pi(x(srt))V_{srt}h\, |\,
\pi(y(sr))V_{sr}k\big)\Delta_G(sr)^{-1}\Delta_H(t)^{-1/2}\rho(sr)^{-1}\notag\\
&\hspace{11cm} dr\, d(sH)\, dt\notag\\
&=\int_H\int_{G/H}\int_H \big(\pi(x(st))V_{st}h\, |\,
\pi(y(sr))V_{sr}k\big)\Delta_G(sr)^{-1}\Delta_H(r^{-1}t)^{-1/2}\rho(sr)^{-1}\notag\\
&\hspace{11cm}
 dr\, d(sH)\, dt. \label{eq-iso2}
\end{align}
Since $\rho(st)\Delta_G(t)=\rho(s)\Delta_H(t)$ for $s\in G$ and
$t\in H$, we have
\[\Delta_G(sr)^{-1}\Delta_H(r^{-1}t)^{-1/2}\rho(sr)^{-1}
=\Delta_G(sr)^{-1/2}\rho(sr)^{-1/2}\Delta_G(st)^{-1/2}\rho(st)^{-1/2},\]
and  an application of
Fubini's Theorem gives
\begin{align*}
\eqref{eq-iso2}
&=\int_{G/H}\Big(\int_H
\pi(x(st))V_{st}h\rho(st)^{-1/2}\Delta_G(st)^{-1/2}\, dt\\
&\hspace{5cm} \Big|\, \int_H
\pi(y(sr))V_{sr}k\rho(sr)^{-1/2}\Delta_G(sr)^{-1/2}\, dr\Big)\,
d(sH),
\end{align*}
which is $\big(\Psi(x\otimes h)\,\big|\,\Psi(y\otimes k)\big)$. Thus $\Psi$ is isometric, and extends to an isometric linear map of $X\otimes_{A\rtimes_\alpha H}\H$ into $L^2(G/H,\H)$.

To see that $\Psi$ has dense range, let   $\xi\in C_c(G/H)$ 
\kaz{be nonzero, and fix $h\in \H$;} 
since $C_c(G/H)\odot\H$ is dense in $L^2(G/H)\otimes\H\cong
L^2(G/H,\H)$, it suffices to approximate $\xi\otimes h:sH\mapsto \xi(sH)h$  in the inductive
limit topology by functions of the form $\sum_{i=1}^n\Psi(x_i\otimes h_i)$ for $x_i\in X_0$ and $h_i\in\H$.  There exists 
\kaz{nonzero}
$f\in C_c(G)$ such that $(\supp f)/H=\supp\xi$ and
\[
\xi(sH)=I(f)(sH):=\int_H f(st)\, dt.
\]
Fix $\epsilon>0$ and a relatively compact open neighbourhood $K$ of $\supp f$. Since $\pi$ is nondegenerate, an approximate identity argument shows that there exists $a\in A$ such that $\|\pi(a)h-h\|<\epsilon/3$, and now an $\epsilon/3$ argument implies that there exists a neighbourhood $O$ of the identity in $G$ such
that
\begin{equation}\label{eq-estimate}
s\in O \Longrightarrow \|\pi(a)V_sh-h\|<\epsilon/\big(\|I(f)\|_\infty \m(\supp f)\big).
\end{equation}
Choose $s_i\in G$ such that $\bigcup_{i=1}^n Os_i$ contains $K$, and use a partition of unity to find functions $\{p_i:1\leq i\leq n\}$ in $C_c(G)$ such that $\supp p_i\subset Os_i\cap K$ and $\sum_i p_i=1$ on $\supp f$. Then $\xi=I(f)=\sum_{i=1}^n I(p_if)$.
Let $x_i=p_if\Delta_G^{1/2}\rho^{1/2}\otimes a$.
Note that $\supp x_i\subset K$, and hence the support of $\sum_{i=1}^n\Psi( x_i\otimes V^*_{s_i}h)$ is contained in $K/H$. For $sH\in K/H$ we have
\begin{align}
\Big\|\xi\otimes  h(sH)-\sum_{i=1}^n&\Psi(x_i\otimes V_{s_i}^*h)(sH)\Big\|\leq\sum_{i=1}^n\int_H |p_i(st)f(st)|\,\|h-\pi(a)V_{sts_i^{-1}}h\|\, dt;\label{eq-estimate2}
\end{align}
the integrand vanishes unless $p_i(st)\neq0$, in which case  $sts_i^{-1}\in O$ and \eqref{eq-estimate} implies that
\[
\eqref{eq-estimate2}
< \frac{\epsilon}{\|I(f)\|_\infty\,\m(\supp f)}\sum_{i=1}^n\int_H |p_i(st)f(st)|\, dt
\leq \epsilon.
\]
We have now proved that $\Psi$ extends to a unitary isomorphism of $X\otimes_{A\rtimes_\alpha H}\H$ onto $L^2(G/H,\H)$.

To finish the proof we verify that $\Psi$ intertwines the two representations.
Let $b$ belong to the 
\kaz{dense}
subalgebra $C_c(G\times G/H,A)$ of $(A\otimes C_0(G/H))\rtimes_{\alpha\otimes\lt} G$. Then
\begin{align*}
\big(((\pi &\otimes
\justMH)\rtimes (V\otimes\lambda^{G/H}))(b)\Psi(x\otimes
h)\big) (sH)\\
&=\int_G(\pi\otimes \justMH)(b(r,\cdot))( V_r\otimes\lambda^{G/H}_r)
(\Psi(x\otimes h))(sH)\, dr\\
&=\int_G \pi(b(r,sH))V_r
\Big(\int_H\pi(x(r^{-1}st))
V_{r^{-1}st}h\Delta_G(r^{-1}st)^{-1/2}\rho(r^{-1}st)^{-1/2}\, dt\Big)\\
&\hspace{10cm}
\rho(r^{-1}s)^{1/2}\rho(s)^{-1/2}\, dr\\
&=\int_G\int_H\pi(b(r,sH))V_r\pi(x(r^{-1}st))
V_{r^{-1}st}h\Delta_G(r^{-1}st)^
{-1/2}\rho(st)^{-1/2}\,dt\, dr \\
&=\int_G\int_H\pi(b(r,sH))\pi(\alpha_r(x(r^{-1}st))) V_{st}h\Delta_G(r^{-1}st)^ {-1/2}\rho(st)^{-1/2}\, dt\, dr\\
&=\int_H \pi\Big( \int_G b(r,sH)\alpha_r(x(r^{-1}st))\Delta_G(r)^{1/2}\, dr
\Big)V_{st}h\Delta_G(st)^{-1/2}\rho(st)^{-1/2}\, dt,
\end{align*}
by Fubini's Theorem.  Upon recognising the inner integral as $b\cdot
x(st)$, where $b\cdot x$ denotes the left action of $b\in(A\otimes C_0(G/H))\rtimes_{\alpha\otimes\lt} G$ on $x\in X$,  this becomes
\[\Psi((b\cdot x)\otimes h)(sH)=\Psi\big((X\dashind(\pi\rtimes U)(b))(x\otimes h)\big)(sH),
\]
as required.
\end{proof}


From now on let $N$ be a closed normal subgroup $N$ of $G$. Our goal is
to decide when a representation of $(A\otimes
C_0(G/N))\rtimes_{\alpha\otimes\lt} G$ is coaction-regular, and then to
combine this with Theorem~\ref{thm-ext-reg} to get 
\kaz{a solution to --- or at least a reformulation of --- }
the original extension problem.

To state our next theorem, we need some notation.  The action $\id\otimes \rt$ of $G/N$ by right translation on $A\otimes C_0(G/N)$ commutes with $\alpha\otimes\lt$, and hence induces an action $\beta$ of $G/N$ on  $(A\otimes C_0(G/N))\rtimes_{\alpha\otimes\lt} G$ such that
\[
\beta_{tN}(b)(s, rN)=b(s, rtN)\ \mbox{ for $b\in C_c(G\times G/N, A)$}.
\]

\begin{thm}\label{thm-fm}
Suppose that $\alpha$ is an action of a locally compact group $G$ on a
$C^*$-algebra~$A$, $N$~is a closed normal subgroup 
\kaz{of~$G$,} 
and $\pi:(A\otimes C_0(G/N))\rtimes_{\alpha\otimes\lt} G\to B(\H)$ is a representation. Then $\pi$
is unitarily equivalent to a coaction-regular representation if and only if there
exists a representation $S$ of $G/N$ on $\H$ such
that $( \pi, S)$ is a covariant representation of $((A\otimes
C_0(G/N))\rtimes_{\alpha\otimes\lt}G, G/N,\beta)$.
\end{thm}

\begin{proof}
If we start with a coaction-regular representation
\[
(\phi\otimes \justMN)\rtimes (W\otimes\lambda^{G/N}):(A\otimes C_0(G/N))\rtimes_{\alpha\otimes\lt} G\to B(\H\otimes L^2(G/N)),
\]
then the right-regular representation $S=1\otimes\nu$ has the required property.  So if $\pi$ is unitarily equivalent to  $(\phi\otimes \justMN)\rtimes (W\otimes\lambda^{G/N})$, we just need to move $1\otimes\nu$ to the space of $\pi$ to obtain $S$.

Conversely, suppose that $(\pi,S)$ is covariant for $((A\otimes  C_0(G/N))\rtimes_{\alpha\otimes\lt}G, G/N,\beta)$, and decompose $\pi$ as $(\pi_1\otimes \mu)\rtimes U$. Notice that $S$ commutes with $U$ and gives a covariant representation $(\pi_1\otimes\mu,S)$ of $(A\otimes C_0(G/N),G/N,\id\otimes\rt)$. We want to find a Hilbert space $\H_2$, an isomorphism $Y:\H\to \H_2\otimes L^2(G/N)$ and a covariant representation $(\pi_2,W)$ of $(A,G,\alpha)$ on $\H_2$ such that
\begin{equation}\label{todo}
Y(\pi_1\otimes \mu,U,S)Y^*=(\pi_2\otimes \justMN,W\otimes \lambda^{G/N},1\otimes\nu).
\end{equation}
(There is a notational subtlety here: the $\pi_1\otimes \mu$ on the left-hand side of \eqref{todo} is the representation of $A\otimes C_0(G/N)=A\otimes_{\max} C_0(G/N)$ associated 
\kaz{to the two}
commuting representations 
\kaz{$\pi_1$ and $\mu$}
on the same space, whereas the $\pi_2\otimes \justMN$ on the right-hand side of \eqref{todo} is the representation of $A\otimes C_0(G/N)=A\otimes_{\sigma} C_0(G/N)$ associated to the representations $\pi_2$ and $\justMN$ on the different spaces $\H_2$ and $L^2(G/N)$.)

Applying the 
\kaz{Stone-von~Neumann}
theorem (as in \cite[Remark~C.35]{tfb}, for example) to the representation $(\mu,S)$ of $(C_0(G/N),G/N,\rt)$ gives a Hilbert space $\H_2$ and an isomorphism $Y:\H\to \H_2\otimes L^2(G/N)$ such that 
\[
Y(\mu,S)Y^*=(1\otimes \justMN, 1\otimes \nu).
\]
For $s\in G$ and $f\in C_0(G/N)$, we compute
\begin{align*}
(1\otimes \lambda_{sN}^*)YU_sY^*(1\otimes M(f))
&=(1\otimes \lambda_{sN}^*)YU_s\mu(f)Y^*\\
&=(1\otimes \lambda_{sN}^*)Y\mu(\lt_s(f))U_sY^*\\
&=(1\otimes \lambda_{sN}^*)(1\otimes M(\lt_s(f)))YU_sY^*\\
&=(1\otimes M(f))(1\otimes \lambda_{sN}^*)YU_sY^*,
\end{align*}
and a similar computation using $U_sS_{tN}=S_{tN}U_s$ shows that $(1\otimes \lambda_{sN}^*)YU_sY^*$ also commutes with every $1\otimes \nu_{tN}$. Thus $(1\otimes \lambda_{sN}^*)YU_sY^*$ belongs to the commutant
\[
\big((1\otimes (\justMN\rtimes\nu))(C_0(G/N)\rtimes G/N)\big)'=\big(1\otimes\K(L^2(G/N))\big)'=B(\H_2)\otimes 1,
\]
and hence there is an operator $W_s$ on $\H_2$ such that $(1\otimes \lambda_{sN}^*)YU_sY^*=W_s\otimes 1$, or equivalently such that $YU_sY^*=W_s\otimes \lambda_{sN}$. Straightforward calculations show that $W$ is a strongly continuous unitary representation of $G$.
Since $Y\pi_1(a)Y^*$ commutes with every $1\otimes M(f)=Y\mu(f)Y^*$ and every $1\otimes \nu_{tN}=YS_{tN}Y^*$, there is also a representation $\pi_2$ of $A$ on $\H_2$ such that $Y\pi_1(a)Y^*=\pi_2(a)\otimes 1$. Now the calculation
\begin{align*}
Y(\pi_1\otimes\mu)(a\otimes f)Y^*&=Y\pi_1(a)\mu(f)Y^*\\
&=(\pi_2(a)\otimes 1)(1\otimes \justMN(f))\\
&=(\pi_2\otimes \justMN)(a\otimes f)
\end{align*}
gives the result.
\end{proof}

\begin{remark}
The argument in the last paragraph of the previous proof is the key step in proving that
\[
\big(A\otimes\K(L^2(G/N)),i_A\otimes \justMN, (i_G\otimes\lambda)\rtimes (1\otimes\nu)\big)
\]
is a crossed product for the system
\[
\big(A\otimes C_0(G/N), G\times (G/N),(\alpha\otimes\lt)\times(\id\otimes\rt)\big).
\]
From this, we deduce that the map
\begin{equation}\label{regrep}
\big((i_A\otimes \justMN)\rtimes(i_G\otimes\lambda)\big)\rtimes (1\otimes\nu)
\end{equation}
is an isomorphism of $\big((A\otimes C_0(G/N))\rtimes_{\alpha\otimes\lt}G\big)\rtimes_\beta (G/N)$ onto $(A\rtimes_\alpha G)\otimes \K(L^2(G/N))$. 
Since the isomorphism \[
(A\otimes C_0(G/N))\rtimes_{\alpha\otimes\lt}G\cong (A\rtimes_\alpha G)\rtimes_{\widehat\alpha|}(G/N)
\]
of \cite[Theorem~A.64]{BE} carries $\beta$ into the dual action $(\widehat\alpha|)^\wedge$ and \eqref{regrep} into the regular representation of $((A\rtimes_\alpha G)\rtimes_{\widehat\alpha|}(G/N))\rtimes_{(\widehat\alpha|)^\wedge} (G/N)$ into 
\[
(A\rtimes_\alpha\K(L^2(G/N))=\K((A\rtimes_\alpha G)\otimes L^2(G/N)),
\]
the injectivity of \eqref{regrep} says precisely that the restriction $\widehat\alpha|$ of the dual coaction is maximal in the sense of \cite{EKQ-mc}. 

The maximality of $\widehat\alpha|$ can also be deduced from the general results in \cite{KQ-fmi}, but of course it is reassuring to have a direct argument. In a previous version of this paper, we deduced from  \cite[Proposition~7.1 and Corollary~7.2]{KQ-fmi} that $\widehat\alpha|$ is maximal, and then proved Theorem~\ref{thm-fm} using Mansfield imprimitivity  for maximal coactions \cite[Theorem~5.3]{KQ-fmi}.
\kaz{We are grateful to the referee of that version for suggesting 
the direct approach used here.}
\end{remark}

Theorems~\ref{thm-ext-reg} and~\ref{thm-fm} gives the following improvement to \cite[Theorem~5]{aHRmackey}. (It is the absence of amenability hypotheses rather than the presence of $A$ which is the real improvement.)

\begin{cor}
\label{thm-main} Suppose that $\alpha$ is an action of a locally compact group $G$ on a $C^*$-algebra $A$
and  $N$ is a closed normal subgroup of
$G$. Let $(\pi,U)$ be covariant representation of $(A,
N,\alpha)$ on $\H$.  Then there exists a representation $V$ of
$G$ on $\H$ such that $(\pi, V)$ is a covariant representation of $(A,G,\alpha)$ and $V|_N=U$ if and only if there exists a
representation $S$ of $G/N$ such that $(X\dashind(\pi\rtimes U), S)$ is
a covariant representation of   $\big((A\otimes
C_0(G/N))\rtimes_{\alpha\otimes\lt} G,G/N,\beta\big)$.
\end{cor}

\end{document}